\documentclass[10pt]{article}
\usepackage{amssymb}
\usepackage{amsfonts}
\usepackage{graphicx}
\usepackage{amsmath}
\usepackage{here}
\usepackage{epsfig}
\usepackage{subfig}
\usepackage{color}
\usepackage{fancybox}
\usepackage{dsfont}
\usepackage[normalem]{ulem}
\usepackage[french,english]{babel}

\definecolor{purple}{RGB}{163,35,142}

\usepackage[comma,numbers,square,sort&compress]{natbib}

\usepackage{multirow}
\usepackage{multicol}

\usepackage[margin=4cm]{geometry}

\usepackage{hyperref}

\usepackage[utf8]{inputenc}

\title{Understanding dynamics of {\it Plasmodium falciparum} gametocytes production: Insights from an age-structured model} 

\author{Ramsès Djidjou-Demasse$^{a,*}$, Arnaud Ducrot$^{b}$, Nicole Mideo$^{c}$, Gaëtan Texier$^{d,e}$ }

\date{	{\small $^{a}$MIVEGEC, Univ. Montpellier, IRD, CNRS, Montpellier, France} \\
	{\small $^{b}$ Normandie Univ., UNIHAVRE, LMAH, FR-CNRS-3335, ISCN, 76600 Le Havre, France} \\
	{\small $^{c}$ Department of Ecology \& Evolutionary Biology, University of Toronto, Toronto, Canada} \\
	{\small $^{d}$ Aix Marseille Univ., IRD, AP-HM, SSA, VITROME, IHU Méditerranée Infection, Marseille, France}\\
	{\small $^{e}$ Centre d'Epidémiologie et de Santé Publique des Armées (CESPA), Marseille, France.} \\
	{\small $^*$Author for correspondence: ramses.djidjoudemasse@ird.fr}}

\begin{document}

\maketitle

\begin{abstract}
Many models of within-host malaria infection dynamics have been formulated since the pioneering work of Anderson {\it et al.} in 1989. Biologically, the goal of these models is to understand what governs the severity of infections, the patterns of infectiousness, and the variation thereof across individual hosts. Mathematically, these models are based on dynamical systems, with standard approaches ranging from $K$-compartments ordinary differential equations (ODEs) to delay differential equations (DDEs), to capture the relatively constant duration of replication and bursting once a parasite infects a host red blood cell. Using malariatherapy data, which offers fine-scale resolution on the dynamics of infection across a number of individual hosts, we compare the fit and robustness of one of these standard approaches ($K$-compartments ODE) with a partial differential equations (PDEs) model, which explicitly tracks the ''age'' of an infected cell. While both models perform quite similarly in terms of goodness-of-fit for suitably chosen $K$, the $K$-compartments ODE model particularly overestimates parasite densities early on in infections when the number of repeated compartments is not large enough. Finally, the $K$-compartments ODE model (for suitably chosen $K$) and the PDE model highlight a strong qualitative connection between the density of transmissible parasite stages ({\it i.e.}, gametocytes) and the density of host-damaging (and asexually-replicating) parasite stages. This finding provides a simple tool for predicting which hosts are most infectious to mosquitoes ---vectors of \emph{Plasmodium} parasites--- which is a crucial component of global efforts to control and eliminate malaria. 

\vspace{0.2in}\noindent \textbf{Key words}. Within-host model, malaria, gametocytemia, parasitemia, infectiousness
\end{abstract}

\section{Introduction}
Malaria has always been a public health problem and, since the discovery of malaria parasites in human blood by Charles Laveran in 1880,  remains so despite more than 100 years of research. Malaria continues to have a significant impact on the world with over 400,000 deaths alone each year \cite{WorldHealthOrganization2019}. It is a vector-borne disease caused by five plasmodial species: {\it Plasmodium falciparum, P. vivax, P. malariae, P. ovale, and P. knowlesi}, with \textit{P. falciparum} being the most pathogenic species infecting humans \cite{KhouryEtAl2018}.

The malaria parasite has a complex life cycle involving sexual reproduction occurring in the insect vector \cite{AlanoCarter1990} and two stages of infection within a human (or animal host), a liver stage \cite{Frevert2004} and blood stage \cite{BannisterMitchell2003}. Human infection starts by the bite of an infected mosquito, which injects the sporozoite form of \textit{Plasmodium} during a blood meal. The sporozoites enter the host peripheral circulation, and rapidly transit to the liver where they infect liver cells (hepatocytes) \cite{Frevert2004}. The parasite replicates within the liver cell before rupturing to release extracellular parasite forms (merozoites), into the host circulation, where they may invade red blood cells (RBCs) to initiate blood stage infection \cite{MillerEtAl2013}. Then follows a series of cycles of replication, rupture, and re-invasion of the RBC. Some asexual parasite forms commit to an alternative developmental pathway and become sexual forms (gametocytes) \cite{RussellEtAl2013}. Gametocytes can be taken up by mosquitoes during a blood meal where they undergo a cycle of sexual development to produce sporozoites \cite{AlanoCarter1990}, which completes the parasite life cycle.

The classical model of within-host parasite multiplication in malaria infections was formulated by Anderson {\it el al.} \cite{AndersonEtAl1989}. This model tracks uninfected red blood cells (RBCs), parasitized RBCs (pRBCs) and merozoites. The pioneer work of Anderson {\it el al.} \cite{AndersonEtAl1989}, has been further developed in several directions including in particular immune response, see for instance \cite{GravenorLloyd1998,Hellriegel1992,HetzelAnderson1996,HoshenEtAl2000,LiEtAl2011,MitchellCarr2010,MolineauxDietz1999,AgustoEtAl2019} for human malaria infection. We also mention discrete-time models such as in \cite{DietzEtAl2006}. Those models use an exponential process to describe the rate of rupture of pRBCs and, as a consequence, then fail to capture realistic lifetimes of the pRBCs on short time scales \cite{Saul1998}. One reason for this is that they are essentially Markovian, {\it i.e.} 'memoryless', a RBC that has been parasitized for 40 hours has the same probability of producing merozoites as {\it e.g.} a RBC parasitized less than a hour ago. Moreover, those models are treating some processes that are likely to be kinda continuous as occurring only in a narrow window ({\it e.g.},  the development of parasites within RBCs and the rupture of pRBC followed by the merozoites release phenomenon). 
% \textbf{\nicole{would it be worth adding a clause here that explains ``, which is important if one wants to study SOMETHING OR OTHER''? Given the history of malaria models, it seems like this simplifications are often `good enough'. So, really emphasizing when that's not likely to be true could be helpful?}}  \textbf{\nicole{Is there another argument here about with a discrete-time model, I guess it would get fairly complicated to track compartments with different durations (Dietz tracks time in steps of t+2...does this make multi-day maturation of gametocytes tricky?) Would another option be to point out that using discrete-time formulations means treating some processes that are likely to be kinda continuous (e.g., immune killing?) as occuring only in a narrow window? I assume this can also be a problem?}} 

To correct this issue, some models of malaria infection include $K$-compartments ordinary differential equations (ODEs) representing a progression through a parasite's developmental cycle, {\it e.g.} \cite{GravenorLloyd1998,SaralambaEtAl2011,ZaloumisEtAl2012,IggidrEtAl2006}, or delay differential equations (DDEs) to capture the time pRBCs take to mature before producing new merozoites, {\it e.g.} \cite{HoshenEtAl2000,KerlinGatton2013,CaoEtAl2019,McKenzieBossert2005,SuEtAl2011}. Other approaches are the use of partial differential equations (PDEs) to track the age-structure of the pRBC population \cite{AntiaEtAl2008,KhouryEtAl2018,CromerEtAl2009,DemasseDucrot2013}. It is shown in \cite{FonsecaVoit2015} that DDEs perform better than the ODEs in representing the dynamics of red blood cells during malaria infection. 

The $K$-compartments ODE model can be interpreted as the application of the method of stages (or the ''linear chain trick'') to the life cycle of pRBC, {\it e.g.} see \cite{IggidrEtAl2006,FengEtAl2007,HurtadoKirosingh2019} and references therein. One problem with the $K$-compartments ODE model is how to decide upon the number of repeated compartments  to capture the realistic dynamics of the pRBCs \cite{GravenorEtAl2002}. Determining the distribution of mean waiting times across compartments for the ODE model is also an issue. If the compartments can be considered equivalent to the developmental stages of pRBCs, then parasites might not spend equal time in each stage. 
% \textbf{\nicole{This is a bit confusing to me. It sounds like the distribution of waiting times within a compartment (but that's just exponential right?) so do you mean the mean waiting time? Or do you mean the distribution of mean waiting times across compartments? if the latter, might be good to justify that too...like, if the compartments can be considered equivalent to developmental stages, then parasites might not spend equal time in each stage (they don't -- I think it's roughly 50\% ring stage, 25\% trophozoite, 25\% schizont)}}

We first introduce both mathematical models (PDE and $K$-compartments ODE) and define the model's parameters and outputs. Next, using gametocyte production as a proxy variable of infectiousness, we compared the model outputs from a PDE stage-structured formulation to those from classical $K$-compartments ODE. Furthermore, the output of both mathematical models is used to qualitatively recover the time course of parasitemia, defined as the proportion of all infected RBCs among the total number of RBCs. Finally, the $K$-compartments ODE model (when $K$ is properly chosen) and the PDE model are used to highlight a strong qualitative connection between gametocyte density and parasitemia.

\section{Material and method}

\subsection{Data and methodology}
Our analysis is based on data collected from malariatherapy taken in \cite{EichnerEtAl2001}. Malaria inoculation was a recommended treatment for neurosyphilis between  1940  and  1963. We also refer to \cite{Chernin1984} for a review paper on malariatherapy and the knowledge gained in the understanding of malaria infection. The data we shall use consist in daily records of gametocyte density for twelve patients. Although malariatherapy has been dismissed for obvious ethical reasons, the advantages to use such data are multiple. Indeed, patients are naive to malaria infection and the dynamics are not perturbed by anti-malarial treatments. Let us notice that such data have been widely used in the literature and in particular to estimate mathematical model parameters. We refer to \cite{EichnerEtAl2001} and the references therein.

The method we shall develop consists of devising a mathematical model to describe the intra-host development of the infection and fitting the model to the available data. The output of the mathematical model will allow us to access various quantities related to the time course of the infection, including parasitemia.

\begin{figure}
\begin{center}
\centerline{\includegraphics[width=1.1\textwidth] {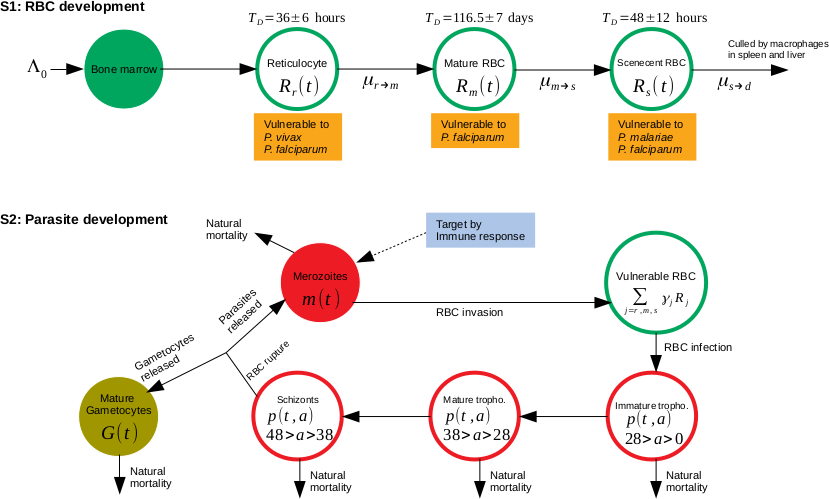}}
\caption{($S_1$) The RBC development chain, ($S_2$) the parasite development chain. $T_D$= average duration ($\pm$ one standard deviation) spent in an RBC age class given in \cite{McQueenEtAl2013}, $\Lambda_0$ is the RBC production rate from the marrow source. In our model, the parameter $1/\mu_{r\to m}$ (resp. $1/\mu_{m\to s }$, $1/\mu_{s\to d}$) is the time spent in RBC reticulocyte (resp. mature, senescent) class.  A continuous parameter $a$ denotes the time since the concerned RBC is parasitized: ring stage ($0<a< 26$ hours), trophozoite ($26<a< 38$ hours) and schizont ($38<a<48$ hours). In the case of \emph{P. falciparum} infection, one has ($\gamma_r =\gamma_m = \gamma_s = 1$) while for \emph{P. vivax} one has ($\gamma_r =1$, $\gamma_m = \gamma_s = 0$) and for \emph{P. malariae}
($\gamma_r = \gamma_m =0$, $\gamma_s = 1$) \cite{PaulEtAl2003}.} \label{flow_diagram}
	\end{center}
\end{figure}

\subsection{Mathematical model}

As discussed above, we now present the mathematical model we shall use to recover parasitemia for twelve patients from observed time courses of gametocyte density. We shall describe the within-host malaria infection coupled with red blood cells (RBCs) production as well as immune effectors. Fig. \ref{flow_diagram} presents the flow diagram of the model considered in this note. Our model is divided into four parts: (i) uninfected RBC (uRBCs) dynamics; (ii) changes in parasite stage or parasite maturity; (iii) Gametocyte production and dynamics and (iv) immune response dynamics.

For uRBCs dynamics, we divide cells into three age classes: reticulocyte (young), mature and senescent. All three ages are vulnerable to \textit{P. falciparum} infection. This can be different for other species of \emph{Plasmodium}. Although we focus in this work on the case of \emph{P. falciparum}, the model described below could be applied to study other species such as  \emph{P. vivax} or \emph{P. malariae}, which have specific RBC-age preferences \cite{PaulEtAl2003}. Such age-structured dynamics for uRBC are well known in the literature, see for instance \cite{McQueenMcKenzie2008}. For the parasites, we consider stage-structured dynamics for their development within pRBC. Here the stage is a continuous variable representing the time since the concerned RBC is parasitized. Such a continuous stage structure will allow us to track the development of parasites within RBCs, but also to have a refined description of the pRBC rupture and of the merozoites release phenomenon. We also emphasize that such a model easily allows for inclusion of anti-malarial treatments acting on only some parasite developmental stages.

\paragraph{Uninfected RBC dynamics.}
We denote by $R_r(t)$, $R_m(t)$ and $R_s(t)$ respectively the density of reticulocytes,
mature RBCs and senescent RBCs at time $t$.

In the absence of malaria parasites, the evolution of circulating
red blood cells is assumed to follow a discrete age maturation system of ordinary differential equations that take the form
\begin{equation}\label{uRBC-model}
\begin{cases}
\frac{dR_r(t)}{dt} = \Lambda_0-\mu_{r\to m}R_r(t),\\
\frac{dR_m(t)}{dt} = \mu_{r\to m}R_r(t)-\mu_{m\to s}R_m(t),\\
\frac{dR_s(t)}{dt} = \mu_{m\to s}R_m(t)-\mu_{s\to d}R_s(t).
\end{cases}
\end{equation}
The parameters $1/\mu_{r\to m}$, $1/\mu_{m\to s}$ and  $1/\mu_{s\to d}$
respectively denote the average duration of RBCs in the reticulocyte, mature and senescent age classes while $\Lambda_0$ represents the normal value of the RBC production from marrow source (i.e. the production rate of RBC). System \eqref{uRBC-model} can also be found in \cite{McQueenMcKenzie2008}.

The parameters of this system are selected from \cite{HetzelAnderson1996,McQueenMcKenzie2008} (see Table \ref{Tab1}) so that in the absence of parasites, the equilibrium age distribution is given by
\begin{equation}\label{eq-initial-RBC}
   \left(R^*_r;R^*_m;R^*_s\right)=\left(62.50; 4853;83.30\right)\times 10^6 \hbox{ cell/ml}. 
\end{equation}
This leads to the homeostatic equilibrium concentration of RBC
$\left(R^*_r+R^*_m+R^*_s\right)$ around $4.99\times10^{9}$ cells/ml which is in the range expected for humans.

\paragraph{Parasite dynamics with stage-structured formulation (PDE model).}
Here we consider the interaction between free merozoites together with the
circulating RBCs. We, respectively, denote by $m(t)$, $p(t,a)$ and $G(t)$ the
density of merozoites, parasitized RBC, and mature gametocytes
at time $t$. The variable $a$ denotes the time since the concerned RBC is parasitized (i.e.
$\int_{a_1}^{a_2}p(t,a)da$ corresponds to the density of pRBC at time $t$ which are infected since the time $a_1<a<a_2$).
The system we shall consider reads as:
\begin{equation}\label{model1}
\begin{cases}
p(t,0)=\beta m(t) \sum\limits_{j=r,m,s}\gamma_jR_j(t),\\
 \partial_t p(t,a)+\partial_a p(t,a)=
-\left(\mu(a)+d_0\right)p(t,a),\\
\dot m(t)=(1-\alpha_G)\int_0^\infty
r\mu(a)p(t,a)da
-\mu_{m}m(t)-\beta m(t)
\sum\limits_{j=r,m,s} \gamma_jR_j(t),\\
\dot G(t) =\alpha_G \int_0^\infty r\mu(a)p(t,a)da -\mu_G
		G(t).
\end{cases}
\end{equation}

We briefly sketch the interpretation of the parameters arising in
\eqref{model1}. Parameters $d_0$, $\mu_m$ and $\mu_G$,
respectively, denote the natural death rates for uRBC, for free
merozoites and for mature gametocytes. Function $\mu(a)$ denotes
the additional death rate of pRBC due to the parasites at stage $a$
and leading to the rupture. The rupture of pRBC at stage $a$ results
in the release of an average number $r$ of merozoites into the
blood stream, so that pRBC then produce, at stage $a$, merozoites at
rate $r \mu(a)$. Together with this description, the quantity
$\int_0^\infty r \mu(a) p(t,a)da$ corresponds to the number of
merozoites produced by pRBC at time $t$. The parameter $\beta$
describes the contact rate between uRBC and free merozoites. Parameters $\gamma_k$ with $k=r,m,s$ describe the age preference of parasites' targets.
Here we shall be concerned in \emph{P. falciparum} infection that does not have any preference for RBC so that
$\gamma_r =\gamma_m = \gamma_s = 1$. However when considering \emph{P. vivax} infection one has
$\gamma_r =1$ and $\gamma_m = \gamma_s = 0$, so that target RBCs mostly consist in reticulocyltes while when \emph{P. malariae} infection is concerned then target RBCs are mostly senescent cells, that is
$\gamma_r = \gamma_m =0$ and $\gamma_s = 1$ \cite{PaulEtAl2003}. The parameter $\alpha_G$ represents the proportion of merozoites from a bursting asexual schizonts that will enter the gametocyte compartment, {\it i.e.}, are ''committed'' to the gametocyte developmental pathway. For simplicity, we have ignored the age structure of gameteocytes and consider $G$ as capturing mature, measurable gametocytes.

\paragraph{Parasite dynamics with $K$-compartments ODE formulation (ODE model).}
For the ODE model formulation, we consider $K$ stages for the pRBC before rupture and set $p=(p_1,p_2,p_3,\cdots,p_K)$, such that $p_j(t)$ denotes the concentration of pRBC at time $t$. Then, setting $\dot z= \frac{\rm{d}z}{\rm{d} t}$ the ODE model writes
\begin{equation}\label{model-ODE}
\begin{cases}
\dot p_1(t)=\beta m(t)\sum\limits_{j=r,m,s}\gamma_jR_j(t) - \left(\mu_1+d_1\right)p_1(t),\\
\dot p_2(t)=\mu_1p_1(t)- \left(\mu_2+d_2\right)p_2(t),\\
\vdots \\
\dot p_K(t)=\mu_{K-1}p_{K-1}(t)- \left(\mu_K+d_K\right)p_K(t),\\
\dot m(t)=(1-\alpha_G)  r\mu_K p_K(t)
-\mu_{m}m(t)-\beta m(t)
\sum\limits_{j=r,m,s} \gamma_jR_j(t),\\
\dot G(t) =\alpha_G  r\mu_K p_K(t) -\mu_G G(t),
\end{cases}
\end{equation}
wherein $1/\mu_i$ the duration of the $i$-stage and $d_i$ the death rate of pRBC. The number of stages $K$ is variable and other parameters and state variables are the same as for the PDE model.

\paragraph{The immune responses.}
Following \cite{DietzEtAl2006}, here we consider two immune responses (IRs) controlling the growth of the parasite population: (i)
an innate IR $S_I(t)$ at time $t$ representing the effect of the pro-inflammatory cytokine cascade and (ii) an adaptive IR $S_A(t)$ at time $t$. The effect of the innate IR is a function of the present parasite (merozoite) density that takes the form
\begin{equation}\label{Sc}
	S_I(t)= \frac{m(t)}{m(t)+S_I^*},
\end{equation}
where $S_I^*$ is the critical parasite density at which the current multiplication factor is reduced by $50\%$.

The adaptive IR is a function of the cumulative parasite density; this function is determined by two host-specific parameters and one constant: (1) $S^*_A$ is the critical
cumulative parasite density at which the current multiplication factor is reduced by $50\%$; (2) $\Delta_0=16$ days is the average delay required by adaptive IR to become effective \cite{DietzEtAl2006}, i.e., for
time $t$ before $\Delta_0$ the cumulative density is set to zero (the adaptive IR has no effect and $S_A(t)=0$ for $t\leq \Delta_0$) and (3)
$\Delta_1=8$ days is the delay that determines the last term in the cumulative density for times $t\ge \Delta_0$,
i.e., 
\begin{equation}\label{Sm}
S_A(t)= \left\{
\begin{split}
&\frac{\int_{\Delta_0}^{t}m(s)ds}{ \int_{\Delta_0}^{t}m(s)ds+S_A^*},\quad   \Delta_0 \le t< \Delta_0+\Delta_1;\\
& \frac{\int_{\Delta_0}^{\Delta_0+\Delta_1}m(s)ds}{ \int_{\Delta_0}^{\Delta_0+\Delta_1}m(s)ds+S_A^*},\quad t \ge \Delta_0+\Delta_1.
\end{split} \right.
\end{equation}

Thus, including these two IR effects, the dynamics of asexual parasite concentration $m(t)$ should be replaced in Models \eqref{model1} and \eqref{model-ODE} respectively by:
\begin{equation}\label{model1-IR}
\dot m(t) = (1-\alpha_G) \int_0^\infty r\mu(a)p(t,a)da - \left(\mu_{m} +\beta 
\sum\limits_{j=r,m,s} \gamma_jR_j(t)+ S_A(t)\right)m(t) -S_I(t),
\end{equation}
and
\begin{equation}\label{model1-IR-ODE}
\dot m(t) = (1-\alpha_G) r\mu_K p_K(t)
- \left(\mu_{m} +\beta 
\sum\limits_{j=r,m,s} \gamma_jR_j(t)+ S_A(t)\right)m(t) -S_I(t).
\end{equation}

\paragraph{Initial conditions.} For both PDE and ODE models, the initial RBCs are assumed to be at their homeostatic equilibrium distribution in the absence of parasites given by \eqref{eq-initial-RBC}, {\it i.e.}, $R_r(0) = R^*_r$; $R_m(0)=R^*_m$; $R_s(0)=R^*_s$. The above models are also assumed to be free of pRBCs at the initial time, and the initial density of malaria parasites is such that $m(0)=m_0$, with $m_0$ a positive constant. These initial conditions are summarized in Table \ref{Tab2}. 

\paragraph{Parasitemia.} The output of both mathematical models can be used to recover the time course of parasitemia, defined  as the proportion of all infected RBC among the total number of RBC.  Using the notation of the model, the parasitemia at time $t$, denoted by $P(t)$ is calculated as follows
\begin{equation}\label{eq-parasitemia}
P(t)=\underbrace{\frac{\int_{0}^\infty p(t,a)da}{\int_{0}^\infty p(t,a)da+\sum\limits_{j=r,m,s} R_j(t)}}_{\text{PDE model}} \text{ or } \underbrace{\frac{\sum_{l=1}^K p_l(t)} {\sum_{l=1}^K p_l(t)+ \sum\limits_{j=r,m,s} R_j(t)}}_{\text{ODE model}}.
\end{equation}

\section{Results}

\subsection{Development of parasites within RBCs and rupture of pRBCs}

An important characteristic of {\it P. falciparum} is the development of parasites within RBCs. The parasite within a RBC then takes an average of 48 hours to mature and release free merozoites. 

With a sequential progression through $K$ stages of parasite maturity before the rupture of the pRBC, the ODE model quantifies the average parasite's development period by 
\begin{equation}\label{eq-duration-sum-mui}
   \frac{1}{\mu_1} + \cdots + \frac{1}{\mu_K}=48 \text{ hours},
\end{equation}
where $1/\mu_i$ is the waiting time across the $i$-stage of maturity. Indeed, the probability of the pRBC of being in the $i$-stage after $a$ hours of infection is then given by $D_i(a)= \mathbb{P}(\tau_i> a)$, where $\tau_i $ denotes the duration within the $i$-th compartment. One assumption of the $K$-compartments ODE model is that variables $\tau_i$,s are independents and exponentially distributed with parameter $\mu_i$ ({\it i.e.} $D_i(a)= e^{-\mu_i a}$, without taking into account other mechanisms such as natural mortality), such that \eqref{eq-duration-sum-mui} is satisfied. Thus, $\sum_{i=1}^K \mathbb{E}[\tau_i] = \sum_{i=1}^K 1/\mu_i =48$ hours.

With the PDE model, the development of parasites within RBCs is characterized by the rupture function $\mu(a)$, which takes the form
$$
\mu(a)=\begin{cases} 0\text{ if $a<48$ hours},\\ \overline{\mu} \text{ if $a\ge 48$ hours}, \end{cases}
$$
where $a$ is the age of the pRBC and $\overline{\mu}$ is a positive parameter. With such formulation, the overall average development period is $\approx 48$ hours as for the ODE model. Indeed, let $D(a)= \exp{\left(- \int_0^a \mu(\sigma) d \sigma\right)}$ the probability that a pRBC remains parasitized after $a$ hours (without taking into account other mechanisms such as the natural mortality). Then, the average parasite's development period is $$\int_0^\infty D(a) da= 48 + \dfrac{1}{\overline{\mu}}.$$
Here we fix, for {\it e.g.} $\overline{\mu}=10$, such that $\int_0^\infty D(a) da= 48 + \dfrac{1}{\overline{\mu}} \approx 48$. The value of $\overline{\mu}$ is therefore not strictly significant as soon as the last approximation holds. 

Consequently, the PDE model formulation allows to continuously track the development of parasites within RBCs and then to have a refined description of the pRBC rupture followed by the merozoites release phenomenon. By contrast, besides the issue of determining the maturation probability $\{D_i\}_{i=1,\cdots,K}$, such a continuous process is quite difficult to capture with the ODE model with $K$ repeated stages.

\begin{figure}
\begin{center}
\centerline{\includegraphics[width=.7\textwidth] {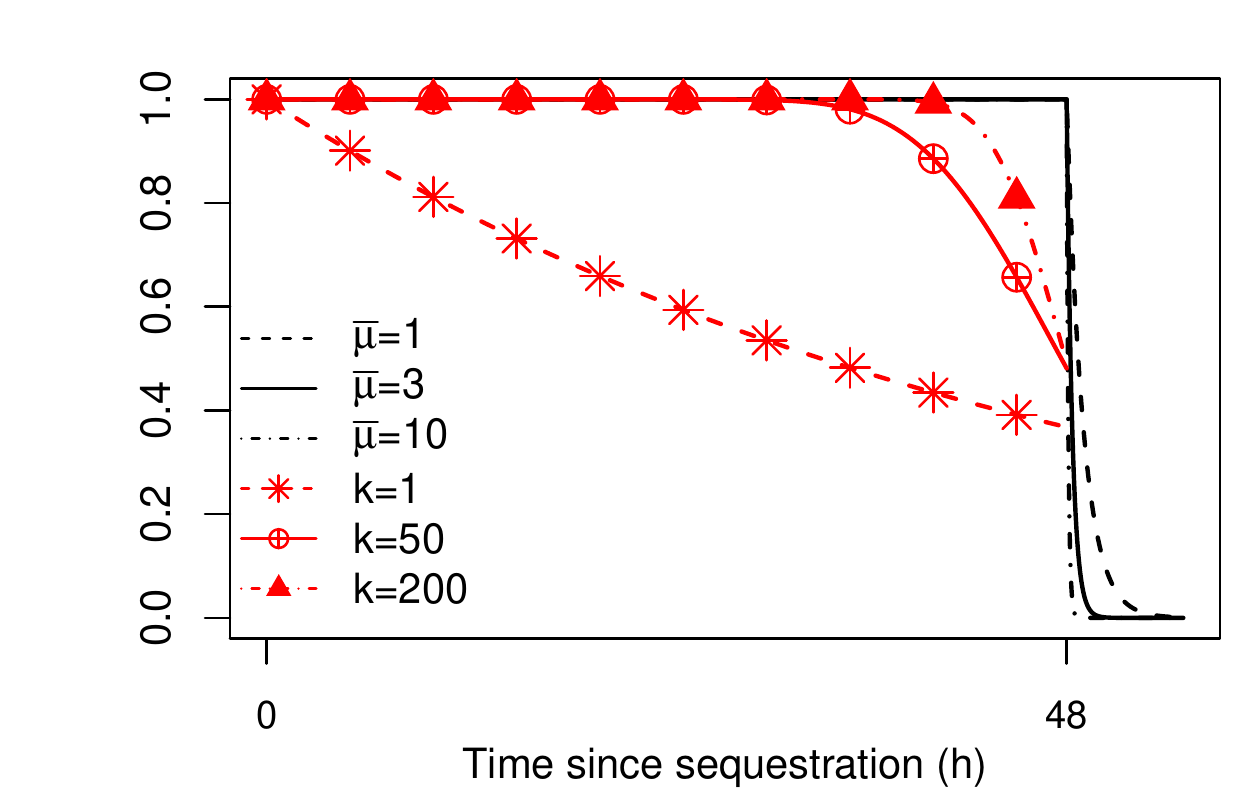}}
\caption{The probability that a pRBC is still parasitized after $a$ hours of infection for both PDE and ODE models. Different values of $K$ correspond to the maturation probability $\mathbb{P}(T_K\ge a)$ for the ODE model, while values of $\overline{\mu}$ are for the function $D(a)$ of the PDE model.} \label{Fig-Sequestration}
\end{center}
\end{figure}

One option for defining the maturation probability $\{D_i\}_{i=1,\cdots,K}$ can be obtained through a ''linear chain trick'' formulation. Indeed, by assuming that the duration of each repeated compartments is the same ({\it i.e.}, $\mu_i=\mu_0$, for all $i=1,\cdots,K$), by \eqref{eq-duration-sum-mui}, we then have $D_i(a)= e^{-a K /48}$, for all $i=1,\cdots,K$. The total duration before rupture becomes $T_K=\sum_{i=1}^K D_i$. Here $D_i$,s are independent and identically distributed with exponential law of parameter $\mu_0=K/48$. Hence, $T_K$ follows a Gamma distribution $\Gamma(K,\mu_0^{-1})=\Gamma\left(K,\frac{48}{K}\right)$. We recover that the mean value of $T_K$ is $48h$ and also that its variance is given by ${\rm var}(T_K)=48^2/K$. The latter quantity tends to $0$ as $K\to \infty$ meaning that $T_K\to 48$ as $K\to \infty$. As a consequence of the above computations, when $K$ is very large the probability that a pRBC remains parasitized after $a$ hours, is approximately given by
$\mathbb{P}(T_K\geq a)\approx 1$ if $a\leq 48$ else $0$, that is close to $D(a)$ when $\overline{\mu}$ is large (Fig. \ref{Fig-Sequestration}). However, the main problem with such formulation is that a very large number of repeated compartments is necessary to capture a quite realistic development process of parasites within a RBC. For instance, the probability that a pRBC is still parasitized after 48 hours of infection is approximately 37\%, 45\%, 48\%, respectively for $K=1,10,50$  (Fig. \ref{Fig-Sequestration}).

% Note that the probability $D(a)= \exp{\left(- \int_0^a \left( \mu(\sigma) +d_0\right) d \sigma\right)}$, that a pRBC remains parasitized after $a$ hours of infection, takes the form $D(a)=e^{-d_0 a}$, for $a\le 48h$ and $D(a)= e^{-d_0 a} \times e^{-(a-48) \overline{\mu}}$, for $a>48h$.
% Therefore, after $48h$ (the average RBC sequestration period), $D$ decreases exponentially fast. Here, we fix $\overline{\mu}=7$, and the value of $\overline{\mu}$ is not strictly significant as soon as the fast decay of $D$ appears after $48h$. 

\subsection{Within-host reproduction number}

The basic reproduction number, usually denoted as $\mathcal{R}_0$, is defined as the total number of parasites arising from one newly pRBC introduced into an uninfected host. It can be used to study the spread of the malaria parasite in an uninfected host, and the parasite will spread if $\mathcal{R}_0>1$. The $\mathcal{R}_0$ of the $K$-compartments ODE and PDE models write 
\begin{equation}\label{R0-ODE-PDE}
\mathcal{R}_0= 
\begin{cases}
\dfrac{ \beta }{\mu_{m} + \beta \sum\limits_{j=r,m,s}  R_j^*}  (1-\alpha_G) r  \prod_{i=1}^K \dfrac{\mu_i}{\mu_i+d_i}  \left( \sum\limits_{j=r,m,s}  R_j^*\right), \text{ ODE} \\
 \dfrac{ \beta }{\mu_{m} + \beta \sum\limits_{j=r,m,s}  R_j^*} (1-\alpha_G) r  \dfrac{\overline{\mu} }{\overline{\mu} +d_0} e^{-48\times d_0} \left( \sum\limits_{j=r,m,s}  R_j^*\right), \text{ PDE}.
\end{cases}
\end{equation}
% \begin{equation}\label{R0-ODE}
%  \text{ODE model} \left \|   \mathcal{R}_0= \dfrac{ \beta }{\mu_{m} + \beta \sum\limits_{j=r,m,s}  R_j^*}  (1-\alpha_G) r  \prod_{i=1}^K \dfrac{\mu_i}{\mu_i+d_i}  \left( \sum\limits_{j=r,m,s}  R_j^*\right), \right.
% \end{equation}
% and 
% \begin{equation}\label{R0-PDE}
%  \text{PDE model} \left \|    \mathcal{R}_0= \dfrac{ \beta }{\mu_{m} + \beta \sum\limits_{j=r,m,s}  R_j^*} (1-\alpha_G) r  \dfrac{\overline{\mu} }{\overline{\mu} +d_0} e^{-48\times d_0} \left( \sum\limits_{j=r,m,s}  R_j^*\right). \right. 
% \end{equation}
We refer to \cite{IggidrEtAl2006,DemasseDucrot2013} for details on the derivation of \eqref{R0-ODE-PDE}.

While the probability for merozoites production for each infection cycle is always $1$ for the ODE model, such probability is $\frac{\overline{\mu} }{\overline{\mu} +d_0}$ for the PDE model. One reason for this is that the $K$-compartments ODE model is essentially Markovian, {\it i.e.} 'memoryless'. With the $K$-compartments ODE model, a RBC that has been parasitized for 40 hours has the same probability of producing merozoites as {\it e.g.} a RBC parasitized less than a hour ago. However, parameter $\overline{\mu}$ can then be chosen such that these probabilities are close to unity. For instance, for the PDE model, $\frac{\overline{\mu} }{\overline{\mu} +d_0} \approx 1$ as soon as $\overline{\mu}$ is sufficiently large compared to $d_0$. Therefore, with the value of $\overline{\mu}=10$ introduced in the previous section, we have $\frac{\overline{\mu} }{\overline{\mu} +d_0} \approx 0.99$.

Finally, one of the main differences between the $\mathcal{R}_0$ expressions of both models is the probability at which pRBCs survive the 48 hours of the parasite's development period for each infection cycle. While such probability is quantified by the term $ e^{-48\times d_0}$  for the PDE model \eqref{R0-ODE-PDE}, it is $\prod_{i=1}^K \frac{\mu_i}{\mu_i+d_i}$ for the ODE model \eqref{R0-ODE-PDE}. However, in some configurations of parameters $\mu_i$,s and for $K$ sufficiently large we can have $\prod_{i=1}^K \frac{\mu_i}{\mu_i+d_i} \approx e^{-48\times d_0}$. For instance, by assuming that the duration of each repeated compartments is the same for the ODE model ({\it i.e.}, $\mu_i=\mu_0$, for all $i=1,\cdots,K$), equality \eqref{eq-duration-sum-mui} gives $K/\mu_0=48$ such that 
\begin{equation*}
\begin{split}
\prod_{i=1}^K \frac{\mu_i}{\mu_i+d_i}= \left(1+d_0/\mu_0\right)^{-K} = \exp\left[-K\ln \left(1+\frac{48d_0}{K}\right)\right] \approx \exp(-48d_0)\text{ if $K\gg 1$}.   
\end{split}
\end{equation*}

\subsection{Fitting the model parameters with data}
The model presented above is solved numerically by using finite volume numerical schemes (implemented with the MatLab Programming Language). The model is then fitted to the data for the time course of gametocytes of the patients. To fit our model, let us observe that most of the parameters are estimated from the literature \cite{AndersonEtAl1989,McQueenMcKenzie2008,EichnerEtAl2001,McQueenEtAl2013,HetzelAnderson1996}. Table \ref{Tab1} provides the values we shall use for the fixed parameters. Three parameters need to be estimated from the data, these are: the proportion of parasitized cells that produce asexual merozoites  ($\alpha_G$), the merozoite initial density ($m_0$), and  the duration of sexual stage ($1/\mu_G$). Additionally to the above three parameters, the number $K$ of compartments is also estimated for the ODE model. These parameters are adjusted from the data for each patient by using a least square method. Basically, we find the values which minimize the difference between the ODE model prediction gametocyte density and the observed data by using MatLab nonlinear least-squares solver {\it lsqcurvefit}. Those optimal parameters for the ODE model are then used to run the PDE model. The superposition of the data and gametocyte density output of the mathematical models are presented in Fig. \ref{figDataModel}, while the estimated parameter values for each patient are given in Table \ref{Table-parameter-fitted}.

\subsection{Comparison of ODE and PDE model outputs}
\begin{figure}
\centerline{\includegraphics[width=\textwidth] {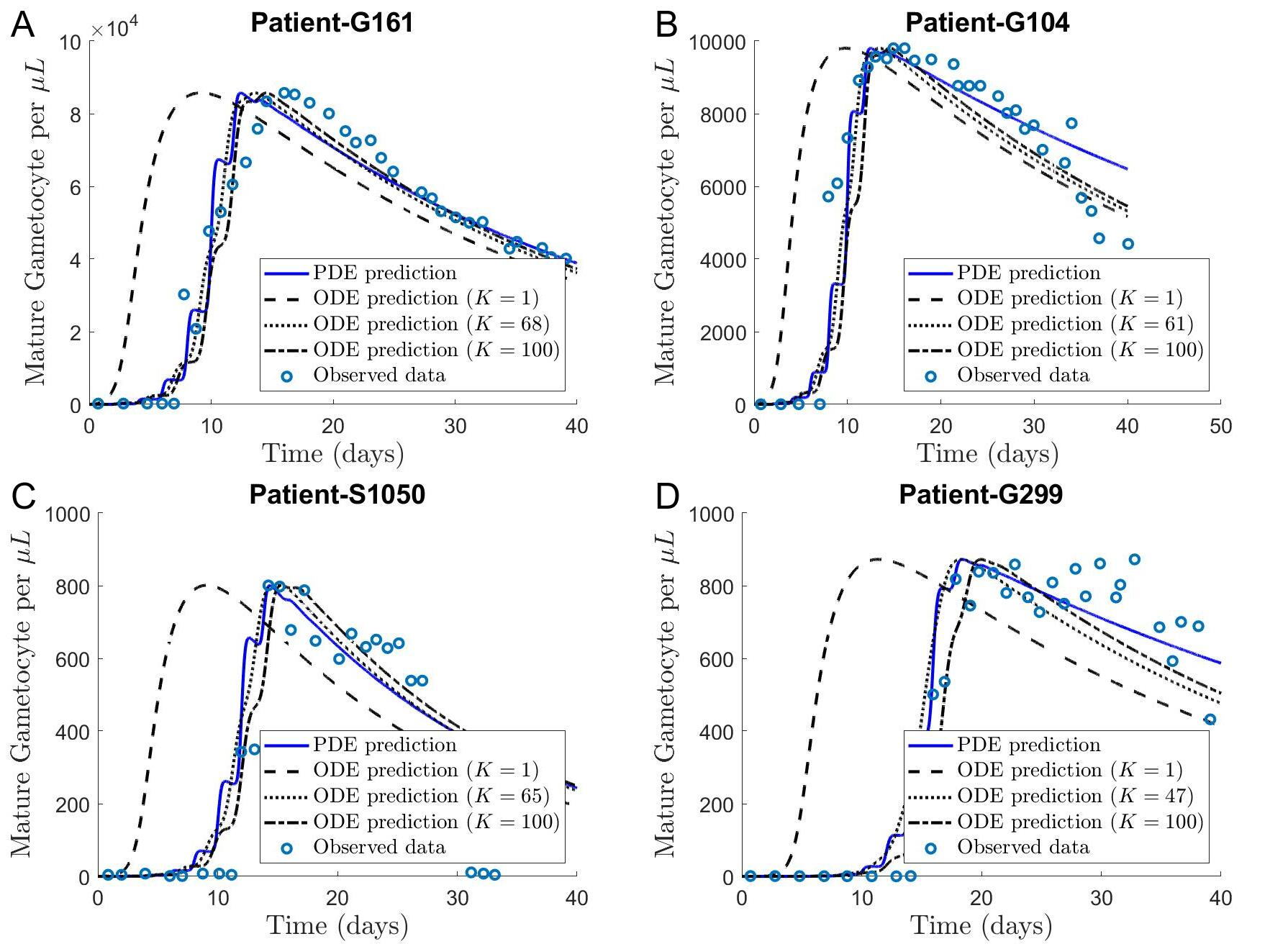}}
% 	\centerline{\includegraphics[width=.4\textwidth] {FigNew/Patient-G161DataGameto} 
% 	\includegraphics[width=.4\textwidth] {FigNew/Patient-G104DataGameto}}
% 	\centerline{\includegraphics[width=.4\textwidth] {FigNew/Patient-S1050DataGameto} \includegraphics[width=.4\textwidth] {FigNew/Patient-G299DataGameto}}
\caption{Comparison with data and mathematical model output for gametocyte density for patients G161, G104, S1050 and G299. The ODE model is illustrated for three values of $K$ ($K=1$, $K=100$, and an intermediate value corresponding to the optimal $K$ for each patient). Comparison for other patients is provided by Figure \ref{figDataModel-SuppMat}.} \label{figDataModel}
\end{figure}

\begin{figure}
\centerline{\includegraphics[width=\textwidth] {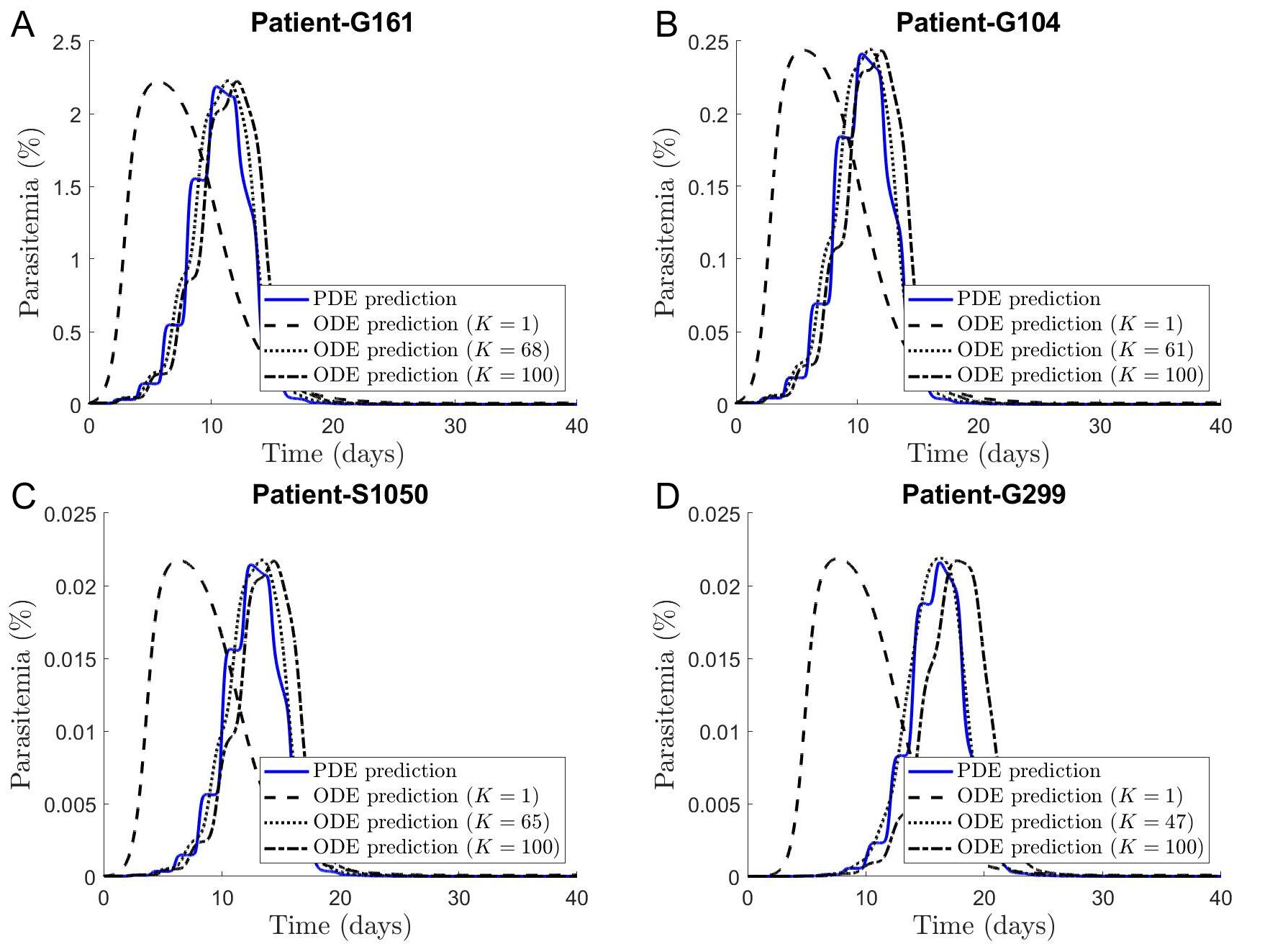}}
% 	\centerline{\includegraphics[width=.4\textwidth] {FigNew/Patient-G161Para} 
% 	\includegraphics[width=.4\textwidth] {FigNew/Patient-G104Para}}
% 	\centerline{\includegraphics[width=.4\textwidth] {FigNew/Patient-S1050Para} \includegraphics[width=.4\textwidth] {FigNew/Patient-G299Para}}
	\caption{The time course of parasitemia (in percentage) for patients G161, G104, S1050 and G299. The ODE model is illustrated for three values of $K$ ($K=1$, $K=100$, and an intermediate value corresponding to the optimal $K$ for each patient). Other patients are provided by Figure \ref{figPara-SuppMat}.} \label{figPara}
\end{figure}

We have presented two modelling frameworks to properly model the within-host infection of malaria. Within this context, we compare a classical model based on ordinary differential equations (ODE) with a model based on partial differential equations (PDE). Our first observation is on the parameterisation of both models. More precisely, a good description of the rupture of pRBC requires at least one additional parameter $K$ for the repeated compartments, see \eqref{model1} versus \eqref{model-ODE}. Such parameter $K$ is necessary to capture the delay in the production or quantification of gametocytes imposed by the development of parasites within RBCs for each infection cycle. This delay in gametocytes production is nicely highlighted by the PDE model formulation (Fig. \ref{figDataModel}). 
% \nicole{\textbf{I think I'm confused by this. It's not sequestration that is being captured here (in Figure 2), it's a delay in the production (or measurement) of gametocytes right? ``Sequestration'' to me means something very specific -- infected RBCs hiding in microvasculature for some part of their life cycle (which it sounds like what is described in the text. But since the data is gametocytes, I don't think it's sequestration that is being captured. I think it would be fine to call this a delay in the production or quantification of gametocytes?}} 
Through a ''linear chain trick'' formulation, it is then possible to find the number $K$ of repeated compartments such that the ODE model can slow down the production dynamics of gametocytes. Indeed, by assuming that the duration of each repeated compartment is the same, we can suitably find the parameter $K$ such that both models perform quite similarly in terms of goodness-of-fit for gametocyte production. Overall, for our dataset, we find between 47 to 68 compartments for a goodness-of-fit of the ODE model (Fig. \ref{figDataModel}). While both models perform quite similarly in terms of goodness-of-fit for a suitable value of $K$, the $K$-compartments ODE model particularly overestimates parasite densities early on in infections when the number of repeated compartments is not large enough (Fig. \ref{figDataModel}). Importantly, the number $K$ of compartments for the goodness-of-fit of the ODE model is quite large and can be variable across individuals (Fig. \ref{figDataModel}). According to the infection dynamics, our comparative results show that the PDE model and the $K$-compartments ODE model (when $K$ is suitably chosen) reproduced, at least qualitatively, the true dynamics of malaria infection parasitemia. Indeed, although we do not have parasitemia real data for patients considered here, but in a qualitative comparison to some studies ({\it e.g.}, \cite{ChildsBuckee2015}), the dynamics of both models seem to mimic qualitatively the parasitemia dynamics when the number of compartments for the ODE model is adequately defined (Fig. \ref{figPara}). 
% \nicole{\textbf{I'm confused about that line too. As i look at figure 3 and I see parasitemia values of like 90\% I worry that cant' be right...I'm I reading this right??? 90\% of RBCs are infected? That is super duper high and I don't think that is what we'd expect for realistic parasitemia values. Have we talked about this before? Have I got stuck on those values before!?}}

%\nicole{\textit{I'm not sure I understand this. Are you talking about infected RBCs sequestering in the microvasculature? This should happen for late stage parasites, but this text makes it sound like they're sequestering for the full 48 hours. I'm not sure how this would explain the 0 estimates for gametocytes observed for ~10 days (if that's what this line is trying to explain). This is a consequence of low infection densities and small fractions of pRBCs making gametocytes right? So, hard to find with early sampling? Happy to chat about this if I'm not making sense!!}} 

\begin{figure}
\centerline{\includegraphics[width=\textwidth] {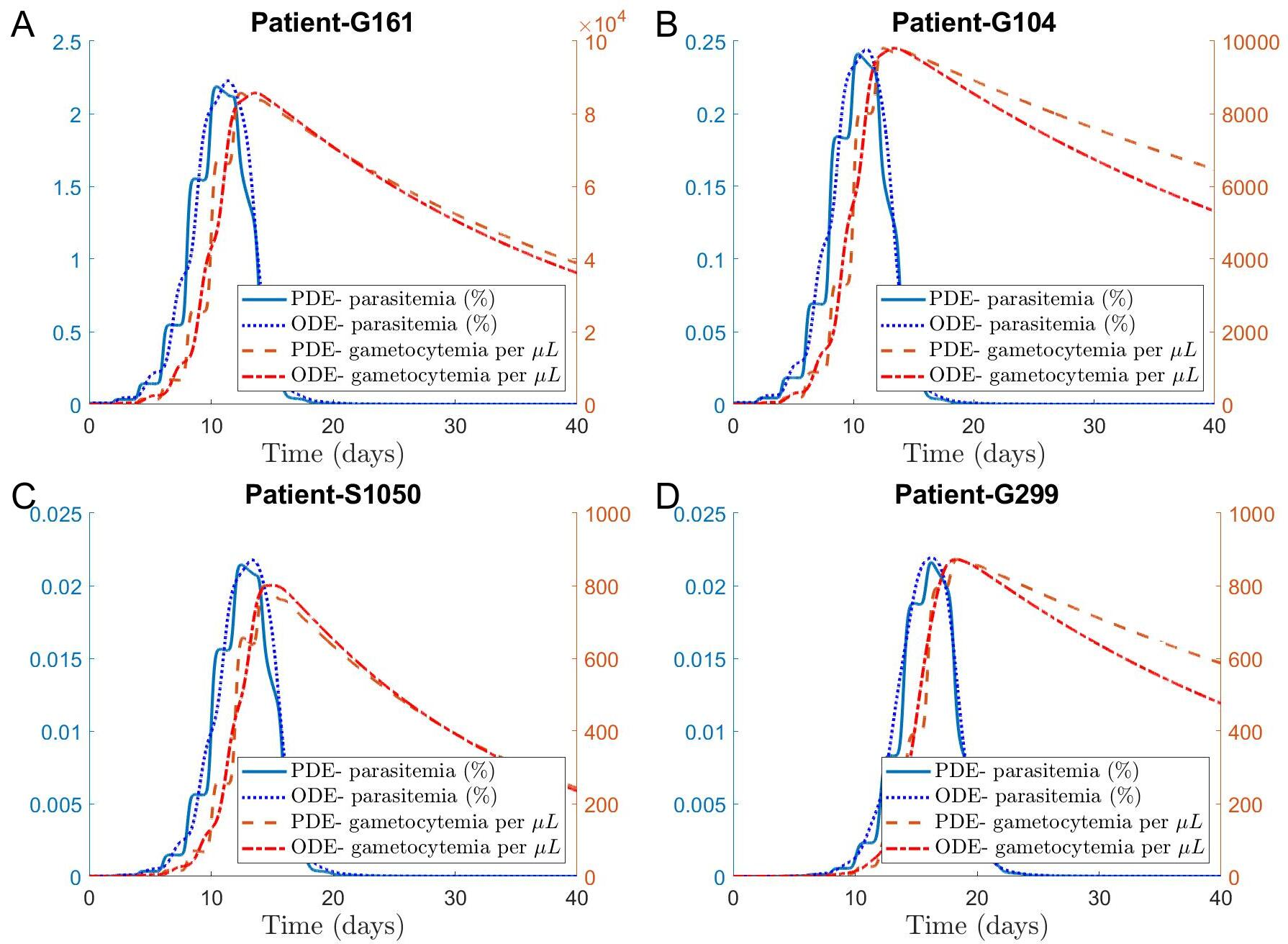}}
% 	\centerline{\includegraphics[width=.4\textwidth] {FigNew/Patient-G161ParaGameto} 
% 	\includegraphics[width=.4\textwidth] {FigNew/Patient-G104ParaGameto}}
% 	\centerline{\includegraphics[width=.4\textwidth] {FigNew/Patient-S1050ParaGameto} \includegraphics[width=.4\textwidth] {FigNew/Patient-G299ParaGameto}}
	\caption{The time course of parasitemia (in percentage) and gametocyte density  computed from the PDE model for patients G161, G104, S1050 and G299. The ODE model is illustrated only for the optimal $K$ of each patient. Other patients are provided by Figure \ref{figGametoPara-SuppMat}.} \label{figGametoPara}
\end{figure}

\subsection{Relationship between parasitemia and gametocyte density}
Our mathematical model has been fitted to the available data for each patient under consideration, which consists of gametocyte densities over time. We now use the output of this mathematical model to recover the time course of parasitemia, defined  as the proportion of all infected RBC among the total number of RBC, see \eqref{eq-parasitemia}. 
% Using the notation of the model, the parasitemia at time $t$, denoted by $P(t)$ can be calculated as follows
% \begin{equation*}
% P(t)=\frac{\int_{0}^\infty p(t,a)da}{\int_{0}^\infty p(t,a)da+\left(R_r(t)+R_m(t)+R_s(t)\right)}.
% \end{equation*}
The time course of parasitemia, $P(t)$, computed from our model are presented in Fig. \ref{figGametoPara} for each patient together with the fitted gametocyte trajectories. As is observed for each patient, the relationship between these curves exhibits two different regimes.
During some period of time $[2,T_0]$, the two curves are increasing with rather similar shape up to a time shift (of length $2$ days). This means that, in this increasing regime, the gametocyte density at time $t$ depends on the parasitemia at time $t-2$, a delay which reflects the life cycle of the parasites inside the RBCs. After this period of increasing parasitemia and gametocyte density, namely after time $T_0$, both curves are decreasing and the shapes seem to depend upon the specific patient considered. To make these comments more quantitative, we introduce the following formula from an estimation of the gametocyte density $G(t)$ from the parasitemia $P(t)$:
\begin{equation}\label{eq-form}
G(t)=\begin{cases} k_1 P(t-2)^{\theta_1}&\text{ if } 2\leq t\leq T_0\text{ days},\\ k_2 P(t)^{\theta_2}&\text{ if $T_0\leq t\leq 30$},\end{cases} 
\end{equation} 
where $k_1$, $k_2$, $\theta_1$ and $T_0$ are four positive parameters while $\theta_2$ is a negative parameter. Let us mention that this 2 days delay between gametocyte density and parasitemia should not to be confused, for example, with the time to distinguish the mature gametocytes via microscopy. Such time is longer than 2 days and can be well captured by the PDE model or the ODE model when $K$ is adequately chosen (Fig. \ref{figDataModel}). While there is a precise biological relationship between parasitemia and gametocytes density at some point in the future, here we seek a robust statistical relationship and so the delay need not match with, {\it e.g.}, gametocyte maturation time.
% \nicole{\textbf{Here should we say something explicit here, like ``While there is a precise biological relationship between parasitemia and gametocyte density at some point in the future, here we seek a robust statistical relationship and so the delay need not match with, e.g., gametocyte maturation time'' (is that convincing?)}}

To determine the unknown parameters $k_1$, $k_2$, $\theta_1$, $\theta_2$ and the changing time $T_0$ for each patient, we perform a least square analysis. More specifically, we adjust these parameters through a logarithmic scale, that is, through the following formula
\begin{equation*}
\log_{10}G(t)=\begin{cases}\log_{10} k_1+\theta_1\log_{10}P(t-2) &\text{ if $2\leq t\leq T_0$},\\ \log_{10} k_2+\theta_2\log_{10}P(t) &\text{ if $T_0< t\leq 30$}.
\end{cases}
\end{equation*} 
To be more precise, our analysis couples the estimate for the time parameter $T_0$, at which the above formula is changing formulation, together with two linear regressions on each part of the graph. We find that parameters $k_1$, $k_2$, $\theta_1$ and $T_0$ are remarkably robust with respect to individuals while $\theta_2$ depends upon each individual. The results of this analysis as well as the estimated parameters are presented in Fig. \ref{figRegression} and summarized by Table \ref{Tab_parameters_fit} for each patient. The quality of the fit is quantified using the coefficient of determination $R^2$ (for linear regression). It is computed for each patient and for the two regimes independently. This adjustment metric is computed using the sample of points induced by the time discretization of the partial differential equation model. For our four cases, this coefficient of determination $R^2$ is approximately $0.99$ in the first part of the curve and even closer to one in the second part of the graph. 

Coming back to the adjusted parameters described in Table \ref{Tab_parameters_fit}, one may observe that the four parameters $k_1$, $k_2$, $\theta_1$ and $T_0$ have robust values with respect to patients while the parameter $\theta_2$ depends on the patient. Using the average values of the adjusted parameters on the set of malariatherapy data, we derive the following clinical formula:
\begin{equation}\label{operational_formula}
G(t)=\begin{cases} 3.843\times 10^7 \cdot P(t-2)^{+1.0304}&\text{ if } 2< t\leq \overline{T}_0\text{ days},\\ 2.981\times 10^9 \cdot P(t)^{-0.0470}&\text{ if } \overline{T}_0\leq t\leq 30\text{ days},\end{cases} 
\end{equation} 
with $\overline{T}_0=14.5636 \pm 0.0064$ days. The relative error for \eqref{operational_formula} is such that $\left|\frac{\Delta G}{G} \right|^2 \le 2.3884\times 10^{-4}+1.0617 \left|\frac{\Delta P}{P} \right|^2$, where $\left|\frac{\Delta P}{P} \right|$ is the relative error on the measurement of parasitemia. Therefore, if $\left|\frac{\Delta P}{P} \right|\le 5\%$, then $\left|\frac{\Delta G}{G} \right| < 5.38\%$. From practical point of view, notice that formula \eqref{operational_formula} can be really useful to estimate the gametocyte density from the parasitemia measurement without necessarily using the quite `complex' mathematical model described in this note.

\section{Discussion}

Many models of within-host malaria infection dynamics have been formulated since the pioneer work of Anderson {\it el al.} \cite{AndersonEtAl1989} in 1989. These models are based on dynamical systems, with standard approaches ranging from ordinary differential equations (ODEs), to delay differential equations (DDEs) or partial differential equations (PDEs). Most ODE model formulations \cite{GravenorLloyd1998,Hellriegel1992,HetzelAnderson1996,HoshenEtAl2000,LiEtAl2011,MitchellCarr2010,MolineauxDietz1999} assume an exponential process to describe the rate of pRBCs rupture and therefore fail to capture realistic lifetimes of the pRBCs. This issue is somewhat corrected when the development of parasites within RBCs and rupture of pRBCs are modeled either by a set of $K$-compartments ODE \cite{GravenorLloyd1998,SaralambaEtAl2011,ZaloumisEtAl2012,IggidrEtAl2006}, or DDEs \cite{HoshenEtAl2000,KerlinGatton2013,FonsecaVoit2015,CaoEtAl2019}. Other approaches are the use of PDEs to track the infection history of a pRBC \cite{AntiaEtAl2008,KhouryEtAl2018,CromerEtAl2009,DemasseDucrot2013}. Using gametocyte production, parasitemia (the proportion of all infected RBC among the total number of RBC) as proxy variables and malariatherapy data, we found that the PDE model and the ODE model perform similarly in terms of goodness-of-fit when a suitable value of $K$ is chosen (Fig. \ref{figDataModel}). Some disadvantages of the ODE model are that the number $K$ of compartments required to achieve a good fit is quite large and can be variable across individuals, and the ODE model particularly overestimates parasite densities early on in infections when the number of repeated compartments is not large enough (Fig. \ref{figDataModel}). Similar comparison holds for the parasitemia dynamics (Fig. \ref{figPara}).

Not least, the $K$-compartments ODE model (for suitably chosen $K$) and the PDE model highlight a strong qualitative connection between gametocyte density and parasitemia. From a practical point of view, such a relation given by \eqref{operational_formula} can be really useful to estimate the gametocyte density from the parasitemia measurement without necessarily using the quite 'complex' mathematical model described in this note.

Here, immune-mediated parasite killing is only considered for merozoites. This choice for immunity targeting merozoites, rather than parasitized red blood cells is mostly because it is a lot easier with our PDE model formulation, particularly in terms of parameterization. However, in some studies, {\it e.g.} \cite{DietzEtAl2006}, parasite levels are not distinguished by merozoite and parasitized red blood cells, such that immunity is acting against merozoite and parasitized red blood cells. Also note that, while there is evidence that the mature live gametocytes evade the immune response clearance pathway,  immune-mediated parasite killing of immature gametocytes has been shown in some studies \cite{BansalKumar2018}.  Finally, taking the best fit parameters and then altering immunity parameters between their minimum, median, and maximum values estimated in \cite{DietzEtAl2006} has very little impact on the model outputs, namely, parasitemia and gametocytes density (figures not shown). However, this can be explained by the fact that (i) merozoites are only short lived and (ii) the relatively short-term validation of the model presented here.

% In some $K$-compartments ODE models ({\it e.g} \cite{IggidrEtAl2006}), the number of parasites produced at each infection cycle is quantified by $\left(r\mu_K p_K\right)$ instead of $\left(\sum_{l=1}^K r\mu_l p_l\right)$ as in models \eqref{model1} and \eqref{model-ODE} for the $\dot m$-equation. While a such formulation seems useful to capture the delay imposed by the parasite sequestration by RBC,  the predicated parasitemia is widely underestimated when $K\ge 2$ (Fig. \ref{fig-only-final-stage}). 

%\paragraph{Predict mosquitoes infection from the parasitemia.}
Reducing infections in mosquitoes---vectors of \emph{Plasmodium} parasites---is a crucial component of global efforts to control and eliminate malaria \cite{AlonsoEtAl2011}. Because a strong correlation exits between the gametocyte density within a host and infectivity of mosquitoes \cite{CollinsJeffery2003,GravesEtAl1988,BousemaDrakeley2011,ChurcherEtAl2013}, progress towards this goal would be bolstered by quantifying gametocytes and identifying highly infectious hosts \cite{StoneEtAl2015,Transmission2017,BousemaDrakeley2017}. %\textit{I commented out the original sentence in case you hate what I did here!}
%Achieving this goal it is necessary to quantified the gametocyteamia in a relevant way, because a strong correlation exits between the gametociteamia and mosquitoes infectivity \cite{CollinsJeffery2003,GravesEtAl1988,BousemaDrakeley2011,ChurcherEtAl2013}.
%However, this is not an easy task and requires more sensitive molecular techniques rather than the conventional microscopy methods \cite{ChurcherEtAl2013,NilssonEtAl2015}. 
On the other hand, parasitemia is easily quantified by light microscopy and therefore is more technically accessible, particularly in regions where malaria is endemic. Therefore, quantifying the relationship between the gametocyte density and parasitemia is of great interest to define more simple tools for the prediction of mosquito infection. The results presented in this note provide one such tool. %\nicole{\textit{I changed that last sentence. Too bold?}}

%\paragraph{Linking the within- and between-host infections level.}
From a public health or population dynamics point of view, the time course of the disease at the between-host level is strongly related to the basic reproduction number (also denoted here by $\mathcal{R}_0$). At the between-host level, the $\mathcal{R}_0$ is defined as the number of secondary infections from a single infected individual introduced in a fully susceptible population. This important metric can be estimated from real data but also using mathematical models. The simplest (deterministic) mathematical model reads as the Ross system of equations from which one can compute this threshold number $\mathcal R_0$ as follows:
\begin{equation*}
\text{Between-host} \left\| \mathcal R_0=\mathcal R_0^{VH}\times \mathcal R_0^{HV}\text{ with }\begin{cases} \mathcal R_0^{VH}=abd_M,\\ \mathcal R_0^{HV}=macd_H.\end{cases} \right.
\end{equation*}
Note that the above $\mathcal R_0$ is for between-host malaria dynamics and is not for the within-host models presented here. In the above formula of $\mathcal R_0$, $m$ represents the number of mosquitoes per person, $a$ denotes the mosquito biting rate, $b$ and $c$ denote the per bite transmission probability respectively from mosquito to human and from human to mosquito, while $d_H$ and $d_M$ correspond respectively to the human recovery and the mosquito death rates. Although more ingredients can been included into the mathematical model, leading to different formulations for $\mathcal R_0$, the above expression contains the main important parameters \cite{Macdonald1955,SmithEtAl2012,MandalEtAl2011,RuanEtAl2008}. Parameters $b$ and $c$ serve as the link between within- and between-host dynamics, since the transmission rates from (to) a host will depend on the dynamics of what is happening within that host (vector). Furthermore, there is a clear relationship between gametocyte density ($G$) and the transmission probability per bite from human to mosquito ($c$) \cite{JohnstonEtAl2013,CollinsJeffery2003,BousemaEtAl2012,BousemaDrakeley2011,ChurcherEtAl2013,StepniewskaEtAl2008,DrakeleyEtAl1999}. From a practical point of view, the parameter $c$ is difficult to estimate in a relevant way. Indeed, an efficient measurement of $c$ requires a good measure of the gametocyte density which is quite difficult to obtain in practice. Indeed, while we can get either gametocytes or parasitemia from microscopy, but with gametocytes tending to be at lower densities (sometimes orders of magnitude lower),  there is a detectability issue.
% \nicole{\textbf{Should we say something here about why? To my mind, you can get either gametocytes or parasitemia from microscopy, so should be equally difficult. But, I suppose, with gametocytes tending to be at lower densities --- sometimes orders of magnitude lower --- there is a detectability issue. If that's what you're thinking too, perhaps it woudl be good to say that, either here or below when you say that parasitemia is easier..}} 
Thus, a simple way to estimate the gametocyte density will help to infer the parameter $c$. Thanks to formula \eqref{eq-form} proposed here, we then have robust relationships between parasitemia (easier to measure) and gametocyte density, at least during the first days of infection (approximately the first two weeks). 
%\nicole{\textit{I'm a bit foggy on the literature here, but I feel like I've read papers that look at the role of asymtpomatic (and possibly low parasitemia) patients as a reservoir of infection for malaria. Lin et al. 2014 Trends in Parasitology is one paper that comes to mind. Anyway, I wonder if it's worth saying something here about the relative role of transmission in later stages of infection, when I imagine parasitemia is lower as immunity is doing its job. It could provide more context for the fact that the relationship is less robust after the first wave.}}

Overall, the proposed model for the dynamics of gametocytes is probably valid for the first asexual wave which last approximately 40 days for each patient in this malariatherapy dataset. This relatively short-term validation is enough for the aim of the current study. However, for the long-term gametocyte dynamics, we need to bring more complexity in the model proposed here. Indeed, the conversion probability of asexual parasites to circulating gametocytes ($\alpha_G$) should be considered to vary among successive waves of asexual parasitaemia to tackle such issue of the long-term gametocyte dynamics \cite{EichnerEtAl2001}, particularly since smear-positive asymptomatic malaria infections detectable by microscopy are an important gametocytes reservoir and often persist for months \cite{LinEtAl2014}.

The robustness of the model proposed here, especially the formula linking parasitaemia and gametocyte density, is only guaranteed during the first two weeks after infection. Beyond this time, this estimate is highly variable from one patient to another. This variability is explained, at least in part, by the variability of the duration of sexual stage ($1/\mu_G$) which governs the decrease in gametocyte density (Table \ref{Table-parameter-fitted}). One interpretation of the variation in this parameter is that there exists variation in how well individuals clear gametocytes or kill them through immune responses \cite{LinEtAl2014,DoolanEtAl2009}. Less variability is expected at the beginning of infection, where the whole system is less constrained by immunity. This is likely to be true for the malariatherapy patient data presented here, since hosts were initially naive. However in high transmission settings, acquired immunity---particularly in older hosts---may obscure the relationship between early gametocyte and asexual parasite densities our work has revealed. Furthermore, fine scale, longitudinal data could assess the applicability of this relationship across settings and age groups, although such data is understandably difficult to obtain.

% ANOTHER COMMENTS ON THE FACT THAT THE MODEL IS VALIDATED ON NAIF PATIENTS ????
% \nicole{\textit{Yeah, good point. Also thinking about what I just typed above. Perhaps in populations where there is more standing immunity, the nice relationship would break down. Also, if there is individual variation in immunity (or other factors), then you might expect the parasites to behave differently in terms of conversion rate (e.g., Greischar et al. 2016 Evolution), which could also change the conversion rate. In any case, you could point out that this relationship is sufficiently straightforward that it would be easy to test its robustness with other data. Though...that may risk a reviewer saying ``if it's so simple, you should do it''. So maybe best to end the whole thing going back to the PDE/ODE comparison --- emphasizing the best way forward for folks modeling malaria. :)}}

Finally, while our work reveals a simple tool for linking aspects of the early dynamics of malaria infections, it also offers specific suggestions for how best to mathematically describe those infection dynamics more broadly. Both PDE and $K$-compartments ODE models have been adopted to capture the subtleties of malaria parasite life cycles in blood-stage infections \cite{KhouryEtAl2018}. Our work provides more understanding among using the PDE model or the ODE model for the within-host human malaria dynamics.
 %results presented in this work clearly highlight the added value of within-host malaria infection of a modeling approach based on PDEs model compared to a classical $K$-compartments ODE.

\begin{table}[!h]
	\begin{center}
		\caption{Fixed model parameters} \label{Tab1} 
		\hspace{-3cm}		\begin{tabular}{llll}
			\hline Parameters & Description (unit) & Values & References \\
			\hline $\Lambda_0$ & Production rate of RBC (RBC/h/ml)& $1.73\times 10^{6}$
			 & \cite{AndersonEtAl1989,McQueenEtAl2013} \\
			$1/\mu_{r\to m}$ & Duration of the RBC reticulocyte stage (h) & $36$ & \cite{McQueenMcKenzie2008} \\
			$1/\mu_{m\to s}$ & Duration of the RBC mature stage (day) & $116.5$
			 & \cite{McQueenMcKenzie2008} \\
			$1/\mu_{s\to d}$ & Duration of the RBC senescent stage (h) & $48$ & \cite{McQueenMcKenzie2008} \\
% 			$1/\mu_{G}$ & Duration of asexual stage & $156$hr& \cite{eichner_genesis_2001,mcqueen_host_2013}\\
			$\beta$ & Infection rate of uRBC (RBC/ml/day) & $6.27\times 10^{-10}$  & \cite{AndersonEtAl1989} \\
			$d_0$ & Natural death rate of uRBC (RBC.day$^{-1}$) & $0.00833$ 
			 & \cite{AndersonEtAl1989} \\
			$\mu_{m}$ & Decay rates of malaria parasites (RBC.day$^{-1}$) & $48$  & \cite{HetzelAnderson1996} \\
			$r$ & Merozoites multiplication factor (dimensionless) & 16   & \cite{AndersonEtAl1989}\\
			$\alpha_G$ & Proportion of sexual merozoites (dimensionless) & 0.05 & \cite{McQueenEtAl2013} \\
			$S^*_I$ & Innate IR density for $50\%$ of parasite killing (cells.$\mu l^{-1}$) & 2,755  & \cite{DietzEtAl2006} \\
			$S^*_A$ & Adaptive IR density for $50\%$ of parasite killing (cells.$\mu l^{-1}$) &
			20.4 &  \cite{DietzEtAl2006} \\
				$\Delta_0$ & Delay required by adaptive IR to be effective (day) &
			$16$ &  \cite{DietzEtAl2006} \\
% 			$\Delta_0$ & Delay required by adaptive IR to be effective & $16$ days& \cite{dietz_mathematical_2006}  \\
			\hline
		\end{tabular}
	\end{center}
\end{table}

\begin{table}[!h]
\begin{center}
\caption{Initial values for the model} \label{Tab2}
\hspace{-1cm}
\begin{tabular}{lll} \hline Variables & Description
& Initial Values \\ \hline
$R_r(0)$ & Population of reticulocytes RBC & $62\times 10^6$ RBC.$ml^{-1}$ \\
$R_m(0)$ & Population of mature RBC & $4.85\times 10^9$ RBC.$ml^{-1}$ \\
$R_s(0)$ & Population of senescent RBC & $83\times 10^6$ RBC.$ml^{-1}$ \\
$p(0,.)$& Population of pRBC for the PDE model &
$0$ cells.$ml^{-1}$ \\
 $p_j$ & Population of pRBC for the ODE màdel &
$0$ cells.$ml^{-1}$ \\
$G(0)$ & Population of mature gametocyte & $0$ cells.$ml^{-1}$ \\
$m(0)$ & Population of malaria parasites & variable \\
			\hline
		\end{tabular}
	\end{center}
\end{table}

\section*{Acknowledgements}
Authors thank Samuel Alizon, Laurence Ayong, Carole Eboumbou and Christophe Rogier for comments and suggestions to improve the manuscript.

\section*{Code availability}
The code (with the MatLab Programming Language) used to simulate the model can be accessed through the Zenodo platform at \url{https://doi.org/10.5281/zenodo.5526271}

\bibliographystyle{abbrv}
\bibliography{MalariaWithinHost.bib}

\newpage

%start pagenumbers
%\clearpage 
%\pagenumbering{arabic} 
\begin{appendix}
%for SUPPLEMENTAL MATERIALS%%%%%%%%%%%%%%
%Table 1 ---> Table S1
\setcounter{table}{0}
\let\oldthetable\thetable
\renewcommand{\thetable}{S\oldthetable}

%Figure 1--> Figure S1
\setcounter{figure}{0}
\let\oldthefigure\thefigure
\renewcommand{\thefigure}{S\oldthefigure}
%%%%%%%%%%%%%%%%%%%%%%%%%%%%%%%%%%%%%%%%%
%{\center{\Huge \textbf {Supporting Information} }}\\ 

%equation
\numberwithin{equation}{section}
\setcounter{equation}{0}

\section{Supplementary figures}

\begin{figure}[H]
\caption{Comparison with data and mathematical model output for gametocyte density for other patients. The ODE model is illustrated for three values of $K$ ($K=1$, $K=100$, and an intermediate value corresponding to the optimal $K$ for each patient).} \label{figDataModel-SuppMat}
\centerline{\includegraphics[width=\textwidth] {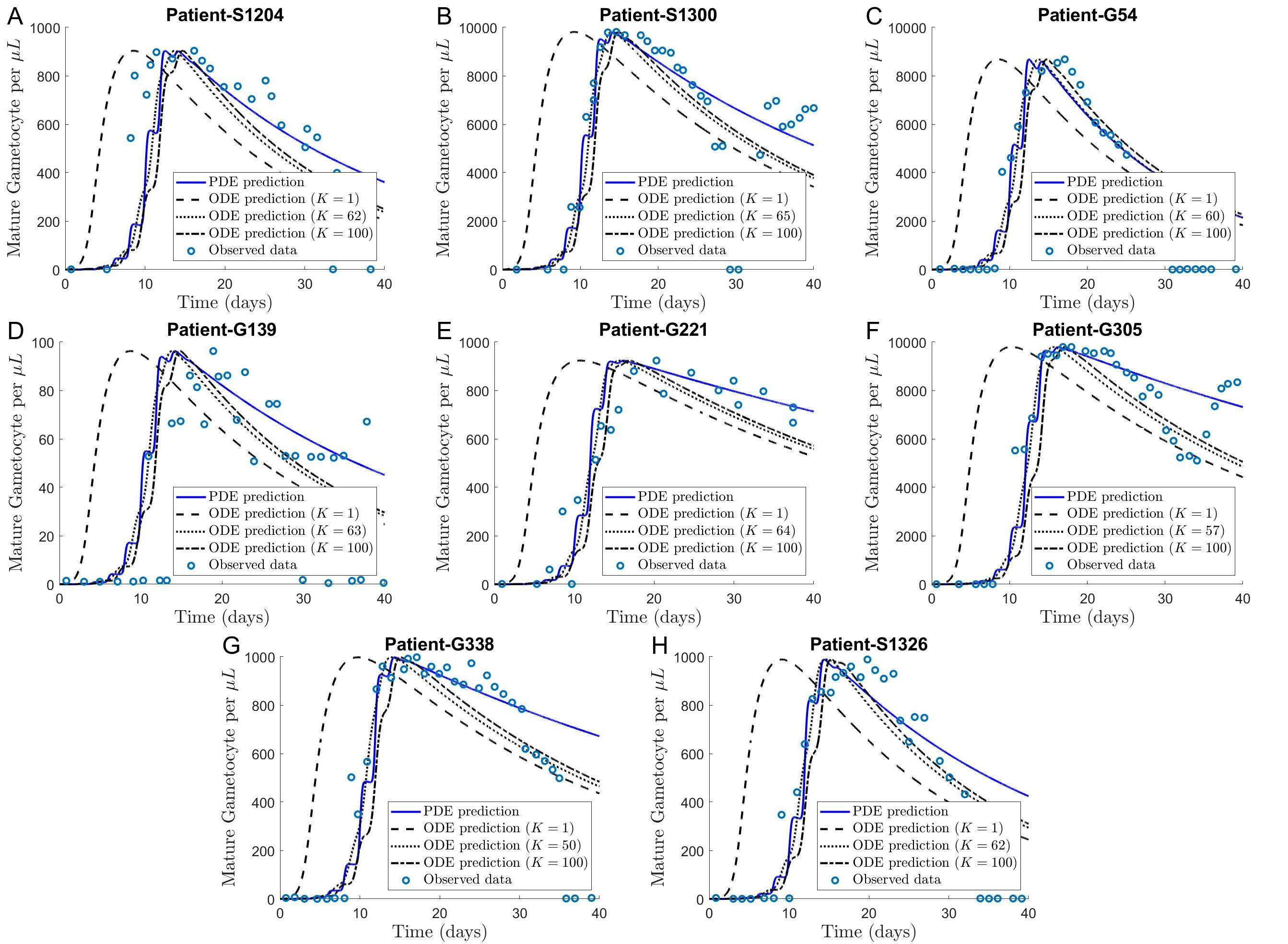}}
% 	\centerline{\includegraphics[width=.3\textwidth] {FigNew/Patient-S1204DataGameto} 
% 	\includegraphics[width=.3\textwidth] {FigNew/Patient-S1300DataGameto} \includegraphics[width=.3\textwidth] {FigNew/Patient-G54DataGameto}}
% 	\centerline{ \includegraphics[width=.3\textwidth] {FigNew/Patient-G139DataGameto} \includegraphics[width=.3\textwidth] {FigNew/Patient-G221DataGameto} \includegraphics[width=.3\textwidth] {FigNew/Patient-G305DataGameto}} 
% 	\centerline{\includegraphics[width=.3\textwidth] {FigNew/Patient-G338DataGameto} \includegraphics[width=.3\textwidth] {FigNew/Patient-S1326DataGameto}} 
\end{figure}

\begin{figure}[H]
\caption{The time evolution of parasitemia (in percentage) density for other patients. The ODE model is illustrated for three values of $K$ ($K=1$, $K=100$, and an intermediate value corresponding to the optimal $K$ for each patient).} \label{figPara-SuppMat}
\centerline{\includegraphics[width=\textwidth] {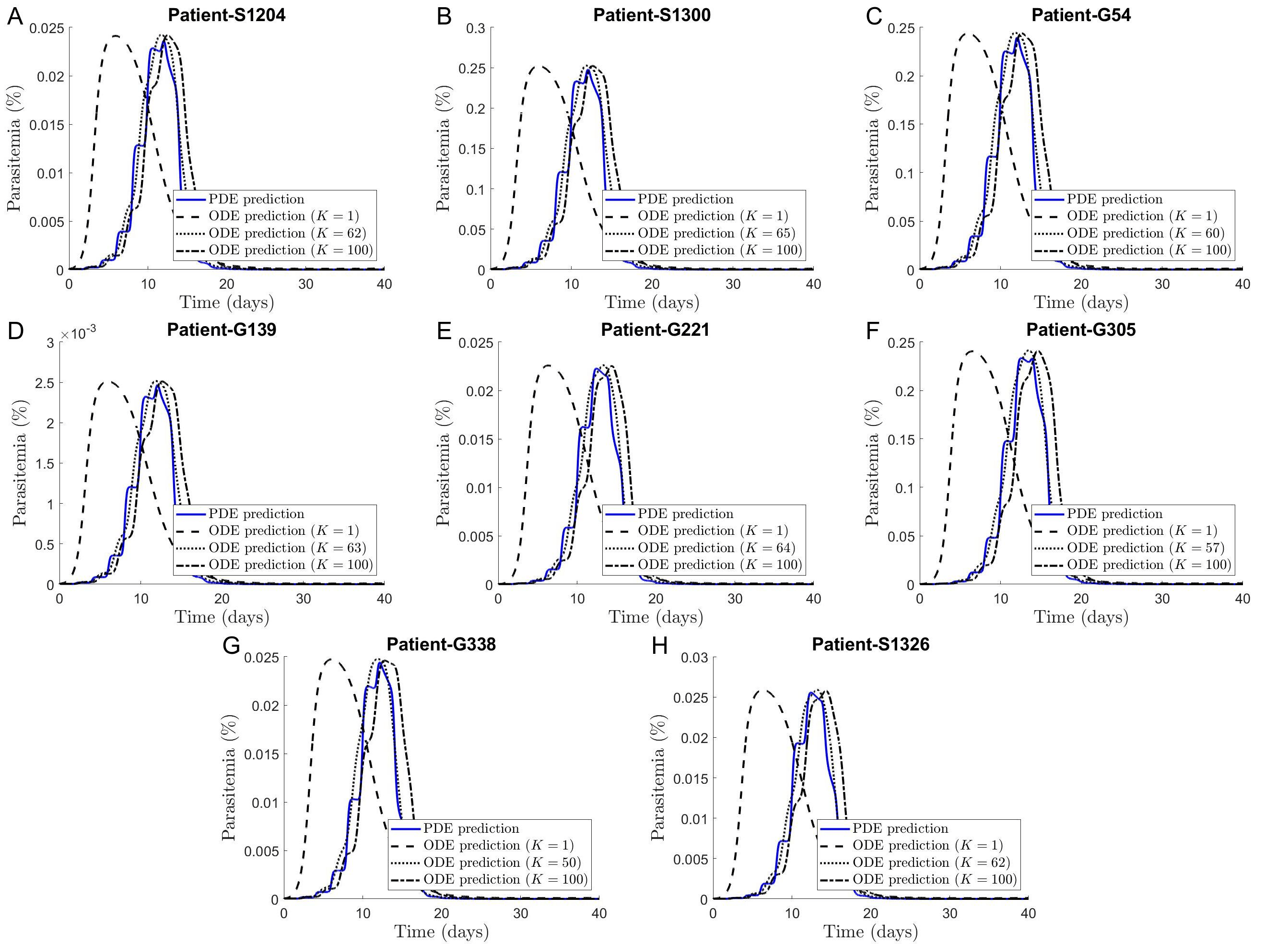}}
% 	\centerline{\includegraphics[width=.3\textwidth] {FigNew/Patient-S1204Para} 
% 	\includegraphics[width=.3\textwidth] {FigNew/Patient-S1300Para} \includegraphics[width=.3\textwidth] {FigNew/Patient-G54Para}}
% 	\centerline{ \includegraphics[width=.3\textwidth] {FigNew/Patient-G139Para} \includegraphics[width=.3\textwidth] {FigNew/Patient-G221Para} \includegraphics[width=.3\textwidth] {FigNew/Patient-G305Para}} 
% 	\centerline{\includegraphics[width=.3\textwidth] {FigNew/Patient-G338Para} \includegraphics[width=.3\textwidth] {FigNew/Patient-S1326Para}} 
\end{figure}

\begin{figure}[H]
\caption{The time evolution of parasitemia (in percentage) and gametocyte density curves computed from the PDE model for other patients. The ODE model is illustrated only for the optimal $K$ of each patient.} \label{figGametoPara-SuppMat}
\centerline{\includegraphics[width=\textwidth] {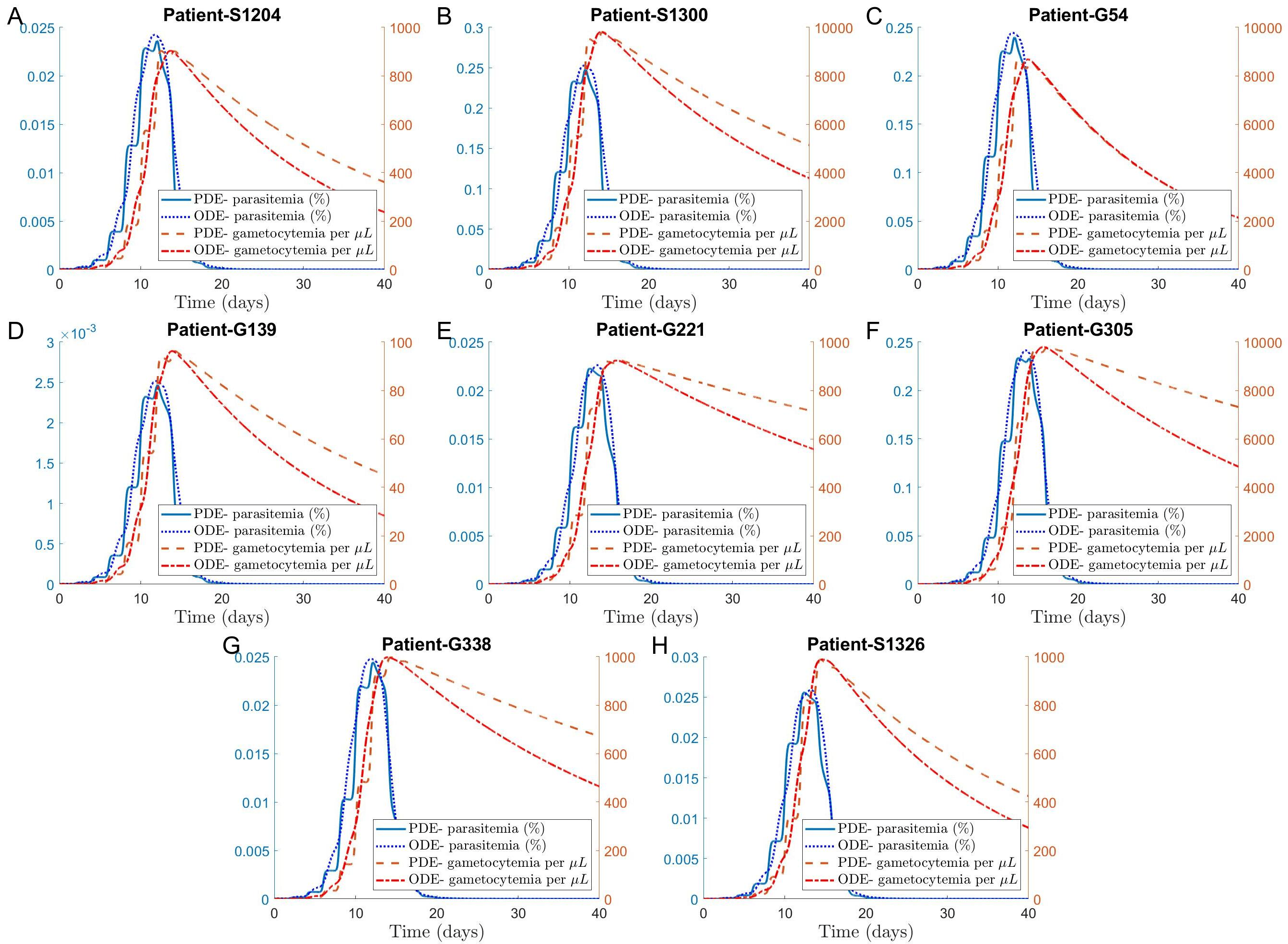}}
% 	\centerline{\includegraphics[width=.3\textwidth] {FigNew/Patient-S1204ParaGameto} 
% 	\includegraphics[width=.3\textwidth] {FigNew/Patient-S1300ParaGameto} \includegraphics[width=.3\textwidth] {FigNew/Patient-G54ParaGameto}}
% 	\centerline{ \includegraphics[width=.3\textwidth] {FigNew/Patient-G139ParaGameto} \includegraphics[width=.3\textwidth] {FigNew/Patient-G221ParaGameto} \includegraphics[width=.3\textwidth] {FigNew/Patient-G305ParaGameto}} \centerline{\includegraphics[width=.3\textwidth] {FigNew/Patient-G338ParaGameto} \includegraphics[width=.3\textwidth] {FigNew/Patient-S1326ParaGameto}} 
\end{figure}

\begin{figure}[H]
\caption{Comparison between model prediction and the linear regression based on formula \eqref{eq-form}} \label{figRegression}
\centerline{\includegraphics[width=\textwidth] {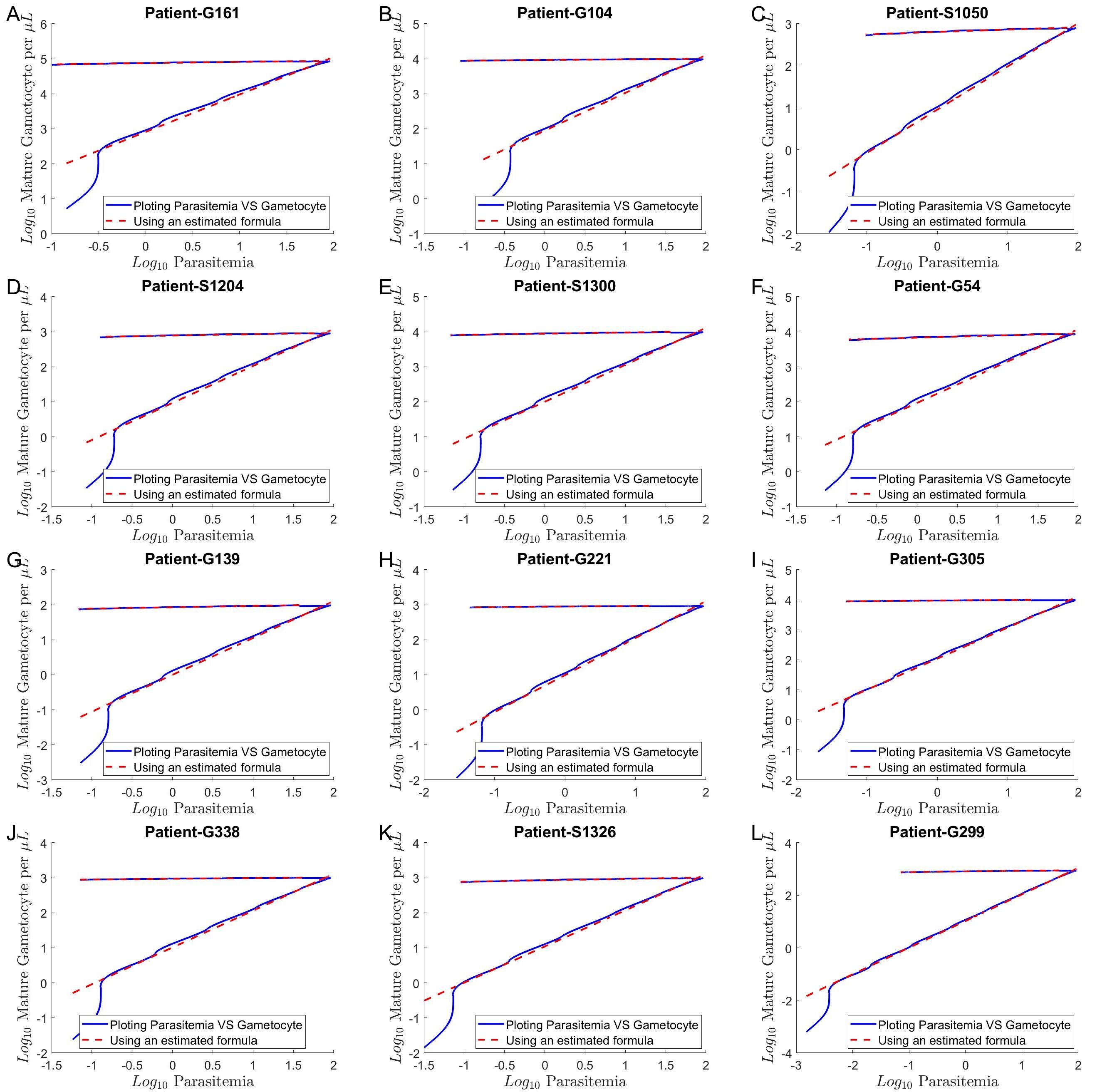}}
%     \centerline{\includegraphics[width=.3\textwidth] {FigNew/Patient-G221Regression} 
% 	\includegraphics[width=.3\textwidth] {FigNew/Patient-G54Regression} \includegraphics[width=.3\textwidth] {FigNew/Patient-S1050Regression}}
% 	\centerline{\includegraphics[width=.3\textwidth] {FigNew/Patient-G299Regression} \includegraphics[width=.3\textwidth] {FigNew/Patient-S1204Regression} 
% 	\includegraphics[width=.3\textwidth] {FigNew/Patient-S1300Regression}}
% 	\centerline{\includegraphics[width=.3\textwidth] {FigNew/Patient-G104Regression} \includegraphics[width=.3\textwidth] {FigNew/Patient-G139Regression} \includegraphics[width=.3\textwidth] {FigNew/Patient-G161Regression}} 
% 	\centerline{\includegraphics[width=.3\textwidth] {FigNew/Patient-G305Regression} \includegraphics[width=.3\textwidth] {FigNew/Patient-G338Regression} \includegraphics[width=.3\textwidth] {FigNew/Patient-S1326Regression}} 
\end{figure}

% \begin{figure}[H]
% \caption{The time evolution of gametocitemia and parasitemia (in percentage) density when the number of parasites produced at each infection cycle is quantified by $\left(r\mu_K p_K\right)$ instead of $\left(\sum_{l=1}^K r\mu_l p_l\right)$.} \label{fig-only-final-stage}
% 	\centerline{\includegraphics[width=.5\textwidth] {FigNew/Patient-G161DataGameto} \includegraphics[width=.5\textwidth] {FigNew/Patient-G161Para}} 
% \end{figure}

\newpage
\section{Supplementary tables}

\begin{table}[!h]
\begin{center}
\caption{Patient-specific parameters estimated from the data %\nicole{\textit{Perhaps you can rearrange this table (and the next one) so that the parameters are columns and the patients are rows? When I look at this table I think that each set of three rows is something really different.}}
}
\label{Table-parameter-fitted}
	\hspace{-1cm}
		\begin{tabular}{c|cccccc} \hline  
			& G54 &G221 &S1050 & S1204 & S1300 & G104 \\\hline 
			$\mu_G (\times 10^{-3})$ &$2.27$ &$0.4632 $ & $1.99$ & $1.5$ & $1.07$ & $0.667$ \\
			$\alpha_G (\times 10^{-8})$ &$13.02$ & $1.206$ & $1.161 $& $1.289$ & $13.47$ & $13.02$  \\
			$m_0 (\times 10^7)$ & $2.5$  & $1$   & $1$  & $2.97$ & $2.5$ & $6$ \\
			\hline \hline 
			& G139 & G161 & G305 &G338 &G299  & S1326 \\\hline 
		$\mu_G (\times 10^{-3})$ & $1.25$ & $1.25$ &$0.526$ &$0.667 $ & $0.80$ & $1.43$  \\
			$\alpha_G (\times 10^{-8})$  & $0.134$& $118.57$ &$12.84$ & $1.321$ & $1.1682$& $1.382$  \\
			$m_0 (\times 10^7)$   & $2.5$  & $5$ & $0.7$  & $2$   & $0.05$  & $1.1$  \\
		\hline
		\end{tabular}
	\end{center}
\end{table}

% \begin{table}[!h]
% \begin{center}
% \caption{Parameters estimated from the data for ODE model }
% \label{Table-parameter-fitted}
% 	\hspace{-1cm}
% 		\begin{tabular}{c|cccccc} \hline  
% 			& G54 &G221 &S1050 & S1204 & S1300 & G104 \\\hline 
% 			$\mu_G (\times 10^{-3})$ &$2.27$ &$46.32 $ & $1.99$ & $1.5$ & $1.07$ & $66.67$ \\
% 			$\alpha_G (\times 10^{-6})$ &$13.02$ & $120.64$ & $116.15 $& $128.89$ & $13.47$ & $13.02$  \\
% 			\hline \hline 
% 			& G139 & G161 & G305 &G338 &G299  & S1326 \\\hline 
% 		$\mu_G (\times 10^{-3})$ & $1.25$ & $1.25$ &$56.63$ &$66.67 $ & $80$ & $1.43$  \\
% 			$\alpha_G (\times 10^{-6})$  & $1341.67$& $1.18$ &$12.84$ & $132.15$ & $116.82$& $138.24$  \\
% 		\hline
% 		\end{tabular}
% 	\end{center}
% \end{table}

\begin{table}[!h]
	\begin{center}
		\caption{Estimated parameters for the shape \eqref{eq-form}} \label{Tab_parameters_fit}
		\hspace{-1cm}
		\begin{tabular}{c|cccc} \hline  
			& G221 &G54 &S1050 & S1204 \\\hline 
			$\log_{10}(k_1)$ &$2.9586\pm 0.0037$ &$2.9326 \pm 0.0039 $ & $2.9520 \pm 0.0037$ & $2.9502 \pm 0.0037$ \\
			$\theta_1$ &$1.0321\pm 0.0028 $ & $1.0312  \pm 0.0030$ & $1.0315  \pm 0.0029 $& $1.0315  \pm 0.0028$  \\
			$\log_{10}(k_2)$ & $2.9437\pm 0.0002 $  & $3.8802 \pm 0.0005$   & $2.9101 \pm 0.0004$  & $2.9148 \pm 0.0003$ \\
			$\theta_2$ &$-0.0141\pm 0.0004$ & $-0.0477 \pm 0.0005$  & $-0.0336 \pm 0.0005$ & $-0.0280  \pm 0.0003$  \\
			$T_0$(days) &$14.7125\pm 0.0075 $ & $13.2917\pm 0.0083 $ &$14.4958\pm 0.0083 $ & $13.4666\pm 0.0025$ \\
			\hline \hline  
		& S1300 &G104 &G139 & G161 \\\hline 
		$\log_{10}(k_1)$ &$2.9861\pm 0.0037$ &$2.9821\pm 0.0041 $ & $2.8024 \pm 0.0037$ & $2.9272 \pm 0.0039$ \\
		$\theta_1$ &$1.0315\pm 0.0028 $ & $1.0304  \pm 0.0031$ & $1.0315  \pm 0.0028 $& $1.0312   \pm 0.0029$  \\
		$\log_{10}(k_2)$ & $3.9442\pm 0.0003 $  & $3.9663 \pm 0.0002$   & $1.9360 \pm 0.0003$  & $4.8774 \pm 0.0002$ \\
		$\theta_2$ &$-0.0336 \pm 0.0003$ & $-0.0170  \pm 0.0002$  & $-0.0336 \pm 0.0003$ & $-0.0214 \pm 0.0003$  \\
		$T_0$(days) &$14.5375\pm 0.0075 $ & $13.2625\pm 0.0083 $ &$14.5625\pm 0.0083 $ & $13.3875\pm 0.0025$ \\
		\hline \hline  
	& G299 &G305 &G338 & S1326 \\\hline 
	$\log_{10}(k_1)$ &$2.9712\pm 0.0034$ &$2.9854 \pm 0.0037 $ & $2.9897 \pm 0.0041$ & $2.9893 \pm 0.0037$ \\
	$\theta_1$ &$1.0301\pm 0.0027 $ & $1.0315  \pm 0.0028$ & $1.0304   \pm 0.0031 $& $1.0315  \pm 0.0028$  \\
	$\log_{10}(k_2)$ & $2.9063\pm 0.0003 $  & $3.9607 \pm 0.0001$   & $2.9739 \pm 0.0002$  & $2.9474 \pm 0.0003$ \\
	$\theta_2$ &$-0.0203\pm 0.0003$ & $-0.0188  \pm 0.0002$  & $-0.0170 \pm 0.0002$ & $-0.0336  \pm 0.0003$  \\
	$T_0$(days) &$20.0083\pm 0.0075 $ & $14.6917\pm 0.0083 $ &$14.2625\pm 0.0083 $ & $14.5458\pm 0.0025$ \\
	\hline
		\end{tabular}
	\end{center}
\end{table}

\end{appendix}

\end{document}